\documentclass[12pt]{amsart}
\newfont{\sheaf}{eusm10 scaled\magstep1}

\usepackage{graphics}
\usepackage{color}

\newcommand{\ra}{\ensuremath{\rightarrow}}

\def\eea{\end{eqnarray*}}
\def\bea{\begin{eqnarray*}}

\newcommand{\Proof}{{\it Proof. }}

\newtheorem{teo}{Theorem}[section]
\newtheorem{df}[teo]{Definition}
\newtheorem{lem}[teo]{Lemma}
\newtheorem{cor}[teo]{Corollary}
\newtheorem{ex}[teo]{Example}
\newtheorem{oss}[teo]{Remark}
\newtheorem{prop}[teo]{Proposition}

\newcommand{\C}{\ensuremath{\mathbb{C}}}

\newcommand{\Z}{\ensuremath{\mathbb{Z}}}
\newcommand{\Q}{\ensuremath{\mathbb{Q}}}
\newcommand{\F}{\ensuremath{\mathbb{F}}}

\newcommand{\N}{\ensuremath{\mathbb{N}}}

\newcommand{\PP}{\ensuremath{\mathbb{P}}}

\newcommand{\FS}{\ensuremath{\mathfrak{S}}}

\title[Absolute Galois acts faithfully on moduli] {The absolute Galois group acts faithfully on the connected
components of the moduli spaces of surfaces of general type}

\author{Ingrid Bauer, Fabrizio Catanese and Fritz Grunewald}

\begin{document}

\date{\today}

\maketitle

\begin{abstract} We show that the Galois group  $Gal(\bar{\Q} /\Q)$ operates faithfully on the set of
connected components of the moduli spaces of surfaces of general type, and also that for each element
$\sigma \in Gal(\bar{\Q} /\Q)$ different from the identity and from complex conjugation, there is a
surface of general type such that
$X$ and the Galois conjugate variety $X^{\sigma}$  have nonisomorphic fundamental groups. The result
was announced by the second author at the Alghero Conference 'Topology of algebraic varieties' in
september 2006. Before the present paper was  actually written, we received a very interesting preprint by
Robert Easton and Ravi Vakil (\cite{e-v}), where it is proven, with a completely different type of examples,
that the Galois group  $Gal(\bar{\Q} /\Q)$ operates faithfully on the set of irreducible components of the
moduli spaces of surfaces of general type.
We also give other simpler examples of surfaces with nonisomorphic fundamental groups
which are Galois conjugate, hence have isomorphic algebraic fundamental groups. 
   \end{abstract}


\section*{Introduction}
In the 60's J. P. Serre showed  in \cite{serre}  that there exists a field automorphism  $\sigma \in
Gal(\bar{\Q} /\Q)
$, and a variety $X$ defined over $\bar{\Q}$ such that 
$X$ and the Galois conjugate variety $X^{\sigma}$  have non isomorphic fundamental groups, in
particular they are not homeomorphic.

In  this note, using completely different methods, we show the following two results, which give a strong
sharpening of the phenomenon discovered by Serre.

\begin{teo} The absolute Galois group $Gal(\bar{\Q} /\Q)$ acts faithfully on the  set of connected
components of the (coarse) moduli spaces of minimal surfaces of general type. 
\end{teo}

\begin{teo} Assume that $\sigma \in Gal(\bar{\Q} /\Q) $ is different from the identity and from complex
conjugation. Then there is a minimal surface of general type $X$ such that
$X$ and $X^{\sigma}$ have non isomorphic fundamental groups.
\end{teo}

In order not to get confused by the above two statements, note that while the absolute Galois group
$Gal(\bar{\Q} /\Q)$ acts on the set  of connected components of the (coarse) moduli spaces of minimal
surfaces of general type, it does not  act on the set of isomorphism classes of fundamental groups of
surfaces of general type. Observe that obviously complex conjugation does not change the isomorphism
class of the fundamental group ($X$ and $\bar{X}$ are diffeomorphic). Now, if we had an action on the set
of isomorphism classes of fundamental groups, then all the normal closure of the  $\Z / 2$ generated by
complex conjugation would act trivially. This would  contradict theorem 2  since complex conjugation does
not lie in the centre of $Gal(\bar{\Q} /\Q)$.

The surfaces we consider in this note are the socalled 'surfaces isogenous to a product' whose weak rigidity
was proven in \cite{isogenous}, and which by definition are the quotient of a product of curves ($C_1
\times C_2$) by the free action of  a finite group $G$.

Therefore our method is strictly related to the socalled theory of  'dessins d' enfants' (see \cite{groth2}).
Dessins d' enfants are,  in view of Riemann's existence theorem (generalized by Grauert and Remmert in
\cite{g-r}),  a combinatorial  way to look at the monodromies of algebraic functions with only three branch
points. It goes without saying that we make here an essential use of Belyi functions (\cite{belyi}) and of
their functoriality.

It would be interesting to obtain similar types of results with examples  involving  rigid varieties belonging
to the class of varieties isogenous to a product. For instance we expect that
 similar results also hold if we restrict
ourselves to consider only the socalled Beauville surfaces (see \cite{isogenous} for the definition of
Beauville surfaces and \cite{cat03},\cite{bcg}, \cite{almeria} for further properties of these).

In the last section we use Beauville surfaces and polynomials with two critical values in order
to produce simple examples of pairs of surfaces with nonisomorphic fundamental groups
which are conjugate under the absolute Galois group (hence the two groups have isomorphic
profinite completions).

\section{Very special hyperelliptic curves}\label{1}

  Fix a positive integer $g \in \N$, $g \geq 3$, and define,
 for any complex number $ a \in \C \setminus \{ -2g,0,1, \ldots , 2g-1 \} $, $C_a$  as the hyperelliptic
curve   of genus $g$ 
$$ w^2 = (z-a) (z + 2g) \Pi_{i=0}^{2g-1} (z-i) $$ branched over 
$\{ -2g, 0,1, \ldots , 2g-1, a \} \subset \PP^1_{\C}$.

\begin{lem} Consider two complex numbers  $a,b$ such that $a \in  \C \setminus \Q $: then $C_a
\cong C_b$ if and only if $a = b$.
\end{lem}

\Proof 

One direction being obvious,  assume that $C_a \cong C_b$. 
Then the two sets with $2g+2$ elements
$B_a : = \{ -2g, 0,1, \ldots , 2g-1, a \}$ and 
$B_b : = \{ -2g, 0,1, \ldots , 2g-1, b \}$ are projectively equivalent over $\mathbb{C}$ (the latter  set $B_b$
has also cardinality $2g+2$ since $C_a \cong C_b$  and $C_a$ smooth implies that also
$C_b$ is smooth).

In fact, this projectivity $\varphi$ is defined over $\Q$, since there are three rational numbers which are
carried into three rational numbers (since $g \geq 2$). 

Since $ a \notin \Q$ it follows that $ b \notin \Q$ and $\varphi$ maps $ B: = \{-2g, 0,1,\ldots 2g-1\}$ to
$ B$, and in particular  $\varphi$ has finite order.  Since $\varphi$ either leaves the cyclical order of
$(-2g, 0,1, \ldots , 2g-1)$ invariant or reverses it,  and $g \geq 3$ we see that there are   $3$ consecutive
integers such that $\varphi$ maps them to $3$ consecutive integers.  Therefore $\varphi$ is either an
integer translation,  or an affine symmetry  of the form $ x \mapsto - x + 2n$. In the former case   $\varphi
= id$, since it has finite order, and in particular, $a = b$. In the latter case it must be $2g + 2n = \varphi
(-2g) = 2g-1$ and $2n = \varphi (0) = 2g-2$,
 a contradiction.

\qed

\begin{oss} The previous lemma  holds more generally under the assumption  that $a,b \in  \C
\setminus
\{ -2g, 0,1, \ldots , 2g-1 \} $, provided $ g \geq 6$.
\end{oss}

\Proof

 The case where $a , b \in \Q$ is similar to the previous one:  $\varphi$ preserves the cyclical order of the
two sets, and we are done if $\varphi (a) = b$ or there are   $3$ consecutive integers which are mapped by
$\varphi$  to $3$ consecutive integers. 

set $B_b$ each consecutive triple of points is a triple of consecutive integers, unless one element in the
triple is $ - 2g$ or $b$. This excludes six triples. Keep  in mind that $ a \in B_a$ and consider all the
consecutive triples of integers in the set $\{0,1,2,3,4,5,6,7, 8 , 9, 10, 11\}$: at most two such triples are not
a consecutive triple of points of $B_a$. We conclude that there is   a triple of consecutive integers in the set
$\{0,1,2,3,4,5,6,7, 8 , 9, 10, 11\}$ mapping to a triple of consecutive integers under
$\varphi$. Then either $\varphi$ is an integer translation $ x \mapsto x + n$, or it is a symmetry  $ x
\mapsto - x + 2 n$. In both case the intervals equal to the convex spans of the sets $B_a$, $B_b$ are sent to
each other by $\varphi$, in particular the  length is preserved and the extremal point are permuted. If  $ a
\in [-2g, 2g-1]$ also  $ b \in [-2g, 2g-1]$ and in the translation case $n=0$, so that $\varphi (x) = x$ and $a
= b$. We see right away that $\varphi$  cannot be a symmetry, because only two points belong to the left
half  of the interval. If  $ a  < -2g $ the interval has length $ 2g-1 -a$, if $ a > 2g-1$ the interval has length $
2g   + a$. We only need to exclude $ a  < -2g $, $ b >  2g-1$, $ 2g-1 -a = 2 g + b$, i.e., $ a = -b-1$. In this case
$b = \varphi (a)$  leads to a contradiction since we should have $b = \varphi (a) = - a + 2n = b + 1 + 2n $.
Else, $\varphi (x ) = x + n$ and  $\varphi (2g-1 ) = b$, $\varphi (2g-2 ) = 2g-1$, hence $n=1$, and
$\varphi (-2g ) = -2 g + 1$, a contradiction.

\qed

We shall assume from now on that $a,b \in \bar{\Q} \setminus \Q$ and that there is a field  automorphism
$\sigma \in Gal(\bar{\Q} /\Q)$ such that $\sigma (a) = b$.  (Obviously, for any  $\sigma \in Gal(\bar{\Q}
/\Q)$ different from the identity, there are  $a,b \in \bar{\Q} \setminus \Q$ with  $\sigma (a) =
b$ and $ a \neq b$.)

\begin{prop} \label{Belyi}
 Let $P$ be the minimal polynomial of $a$ and consider the field $L:=
\mathbb{Q}[x]/(P)$. Let $C_x$ be the hyperelliptic curve over $L$  
$$ w^2 = (z-a) (z + 2g) \Pi_{i=0}^{2g-1} (z-i) .$$ Then there is a rational function
$F_x:C_x
\rightarrow
\PP^1_L$ such that for each $a \in \C$ with $P(a) = 0$ it holds that  the rational function $F_a$ (obtained
under the specialization $ x \mapsto a$) is a Belyi function
 for $C_a$.
\end{prop}

\Proof Let $f_x : C_x \rightarrow \mathbb{P}^1_L$ be the hyperelliptic involution, branched in $\{-2g, 0,1,
\ldots , 2g-1, x \}$. Then $P \circ f_x$ has as critical values:

\begin{itemize}
\item the images of the critical values of $f_x$ under $P$, which are $\in \mathbb{Q}$,
\item the critical values $y$ of $P$, i.e. the zeroes of the discriminant $h_1(y)$ of $P(z) - y $ with respect 
to the variable $z$.
\end{itemize}

$h_1$ has degree $deg(P) - 1$, whence, inductively, we obtain $\tilde{f}_x := h \circ P \circ f$ whose
critical values are all contained in $\Q \cup \{ \infty \}$. If we take any root $a$ of $P$, then obviously
$\tilde{f}_a$ has the same critical values.

Let now $r_1, \ldots , r_n \in \Q$ be the (pairwise distinct) finite critical values of
$\tilde{f}_x$. We set:
$$ y_i := \frac{1}{\Pi_{j \neq i} (r_i - r_j)} \ .
$$

Let $N \in \N$ be a positive integer such that $m_i:=Ny_i \in \mathbb{Z}$. Then we have that the rational
function
$$ g(t):= \Pi_i (t-r_i)^{m_i} \in \mathbb{Q}(t)
$$ is ramified at most in $\infty$ and $r_1, \ldots r_n$. In fact, $g'(t)$ vanishes exactly when the
logarithmic derivative $ G(t) : = \frac{g'(t)}{g(t)} = \sum_i m_i (\frac{1}{t-r_i})$ has a zero or a pole, but
the choice made yields a zero of order $n$ at $\infty$. 

Therefore the critical values of $g \circ \tilde{f}_x$ are at most $0, ~ \infty, ~ g(\infty)$ (for details see
\cite{wolfart}).

We set $F_x : = \Phi \circ g \circ \tilde{f_x}$ where $\Phi$ is the affine map $ z \mapsto 
 g(\infty)^{-1} z$, so that the critical valus of $F_x$ are equal to $ \{ 0,1,\infty\}$.
 It is obvious by our construction that for any root
$a$ of $P$,
$F_a$ has the same critical values as $F_x$, in particular, $F_a$ is a Belyi function for $C_a$.

\qed

Since in the sequel we shall consider the normal closure  (we prefer here, to avoid confusion, not to use the
term 'Galois closure' for the geometric setting) $ \psi_a : D_a \to \PP^1_{\C}$  of each of the functions
$F_a : C_a \to \PP^1_{\C}$, we recall in the next section the 'scheme theoretic' construction of the normal
closure.

\section{Effective construction of normal closures}\label{nc}

In this section we consider algebraic varieties  over the complex numbers, endowed with their Hausdorff
topology, and, more generally,
 'good' covering spaces (i.e., between topological spaces which are arcwise  connected and semilocally
simply connected).

\begin{lem} Let $\pi : X \rightarrow Y$ be a finite 'good'  unramified covering space of degree $d$ 
between connected spaces $X$ and $Y$. Then the normal closure $Z$ of $\pi : X \rightarrow Y$ is
isomorphic to any connected component of 
$$W : = W_{\pi}:= (X \times _Y \ldots \times _Y X) \backslash \Delta \subset X^d
\backslash \Delta ,$$ where $\Delta := \{(x_1, \ldots, x_d) \in X \times _Y \ldots \times _Y X | 
 \exists  i\neq j ~, \   x_i =  x_j\}$ is the big diagonal.
\end{lem}

\Proof  Choose base points $x_0 \in X$, $y_0 \in Y$ such that $\pi (x_0) = y_0$ and denote by $F_0 $ the
fibre over $y_0$, $F_0 : = \pi^{-1}(\{y_0 \})$.

We consider the monodromy $\mu : \pi_1(Y, y_0)
\rightarrow
\mathfrak{S}_d = \mathfrak{S}( F_0)$ of the unramified covering $\pi$. The monodromy of $\phi: W
\ra Y$ is induced by the diagonal product monodromy $\mu^d : \pi_1(Y, y_0) \ra
\mathfrak{S}( F_0^d)$,  such that, for  $(x_1, \ldots , x_d) \in F_0^d$, we have $\mu^d(\gamma) (x_1,
\ldots , x_d) = (\mu(\gamma)(x_1), \ldots ,
\mu(\gamma)(x_d))$.

It follows that the monodromy of $\phi: W
\ra Y$,  $  \mu_W : \pi_1(Y, y_0) \ra
\mathfrak{S}( \mathfrak{S}_d )$ is given  by left translation
$\ \mu_W (\gamma) (\tau) = \mu (\gamma)\circ (\tau) $.

If we denote by $G := \mu ( \pi_1(Y, y_0) )\subset \mathfrak{S}_d $ the monodromy group, it follows right
away that  the components of $W$ correspond  to the cosets $G \tau$ of  $G$. Thus all the components
yield   isomorphic covering spaces.

\qed

The theorem of Grauert and Remmert (\cite{g-r}) allows to extend the above construction to yield normal
closures of morphisms between normal algebraic varieties.

\begin{cor}\label{ncram} Let $\pi : X \rightarrow Y$ be a finite map between normal projective varieties, let
$B \subset Y$ be the branch locus of $\pi$ and set $X^0 := X \setminus \pi^{-1}(B)$,
$Y^0 := Y \setminus B$. 

If $X$  is connected, then the normal closure $Z$ of $\pi$ is isomorphic to any connected component of the
closure of $W^0 := (X^0 \times _{Y^0} \ldots \times _{Y^0} X^0) \backslash \Delta$ in the normalization
$W^n$  of $ W : = \overline{(X \times _Y \ldots \times _Y X) \setminus \Delta}$.

\end{cor}

\Proof

The irreducible components of $W$ correspond to the connected components of $W^0$, as well as to the
connected components $Z$ of $W^n$.  So, our component $Z$ is the closure of a connected component
$Z^0$ of $W^0$. We know that the monodromy group $G$ acts on $Z^0$ as a group of covering
transformations and simply transitively on the fibre of $Z^0$ over $y_0$: by normality the action  extends
biholomorphically to $Z$, and clearly $ Z / G \cong Y$.

\qed

\section{Connected components of moduli spaces associated to very special hyperelliptic curves}

Let $a$ be an algebraic number, $g \geq 2$, and consider as in section \ref{1} the hyperelliptic curve
$C_a$ of genus $g$ defined by the equation
$$ w^2 = (z-a) (z + 2g) \Pi_{i=0}^{2g-1} (z-i) .$$  Let $F_a : C_a \rightarrow \PP^1$ be the Belyi function
constructed in proposition \ref{Belyi} and denote by $\psi_a :D_a \rightarrow \PP^1$ the normal closure
of $C_a$ as in corollary \ref{ncram}.

\begin{oss} We denote by $G_a$ the monodromy  group of $D_a$ and observe that there is a subgroup
$H_a \subset G_a$ acting on $D_a$ such that $D_a /H_a \cong C_a$.
\end{oss} 

 We choose a monodromy representation
$\mu : \pi_1(\PP^1 \setminus \{0, 1, \infty \}) \rightarrow G_a$  corresponding to the normal ramified
covering $\psi_a : D_a \rightarrow \PP^1$ and we denote the images of geometric loops around $0, ~ 1, ~
\infty$ by $ \tau_0, ~ \tau_1, ~ \tau_{\infty}$. Then we have that
$G_a$ is  generated by $ \tau_0, ~ \tau_1, ~ \tau_{\infty}$ and
$\tau_0 \cdot \tau_1 \cdot \tau_{\infty}= 1$.

Fix now another integer $g' \geq 2$ and consider all the possible smooth complex curves
$C'$ of genus $g' $. Observe that the fundamental group of $C'$ is  isomorphic to the standard group 
$$\Pi_{g'} : = \langle \alpha_1,
\beta_1, \ldots ,
\alpha_{g'}, \beta_{g'} | \Pi_{i=1}^{g'} [ \alpha_{i}, \beta_{i}] = 1 \rangle.$$

Since $g' \geq 2$ there are epimorphisms  (surjective homomorphisms) $ \rho \colon \Pi_{g'} \ra  G_a$.
For instance it suffices to consider the epimorphism 
$\theta : \Pi_{g'} \ra \F_{g'}$ from $\Pi_{g'}$ to the free group $\F_{g'}:=<\lambda_1, \ldots ,
\lambda_{g'}>$ in $g'$ letters given by $\theta(\alpha_i) = \theta(\beta_i) = \lambda_i$, $\forall$ $1
\leq i \leq g'$, and to 
 compose $\theta$ with the surjection $\phi : \F_{g'} \rightarrow G_a$, given by $\phi (\lambda_1) =
\tau_0$, $\phi(\lambda_2) = \tau_1$, and $\phi (\lambda_i)= 1$ for $3\leq i \leq g'$.

Consider all the possible  epimorphisms   $ \rho \colon \Pi_{g'} \ra  G_a$. Each such $\rho$ gives a
normal unramified covering $D' \rightarrow C'$ with  monodromy   group $G_a$.

\begin{df} Let $\frak M_a$ be the subset of the moduli space of surfaces of general type given by surfaces
$S \cong (D_a \times D' )/G_a$, where $D_a, D'$ are as above and the group $G_a$ acts by a diagonal
action.
\end{df}

From \cite{isogenous} and especially Theorem 3.3 of \cite{cat03} it follows:

\begin{prop} For each $ a \in \bar{\Q}$,  $\frak M_a$ is a union of connected components of the moduli
spaces of surfaces of general type.

Moreover, for $\sigma \in Gal ( \bar{\Q} / \Q)$, $\sigma (\frak M_a )= \frak M_{\sigma(a)}$.
\end{prop}

\Proof Since $D_a$ is a triangle curve, the pair $ (D_a, G_a)$ is rigid, whereas, varying $C'$ and
$\rho$, we obtain the full union of the moduli spaces for the pairs $(D', G_a)$, corresponding to the
possible free topological actions of the group $G_a$ on a curve $D'$ of genus $ |G_a| (g'-1) + 1$.

Thus, the surfaces $S \cong (D_a \times D' )/G_a$  give, according to the cited theorem 
  3.3 of \cite{cat03},  a union of connected components of the moduli spaces of surfaces of general type.

Choose now  a surface  $S$ as above (thus, $[S] \in \frak M_a$) and apply the field automorphism
$\sigma \in Gal ( \bar{\Q} / \Q)$ to a point of the Hilbert scheme corresponding to the 5-canonical image
of $S$ (which is isomorphic to $S$, since the canonical divisor of $S$ is ample). We obtain a surface which
we denote by $  S ^{\sigma}$.

By taking the fibre product of $\sigma$ with  $D_a \times D' \ra S$ it follows that 
$ S ^{\sigma}$ has an \'etale covering with group $G_a$ which is the product
$ (D_a )^{\sigma} \times ( D') ^{\sigma} $.

Since $(C_a)^{\sigma} = C_{\sigma(a)} $ and since $\sigma (a)$ corresponds to another embedding of the
field $L$ into $\C$, it follows that $(F_a) ^{\sigma} = F_{\sigma(a)} $, whence $(D_a )^{\sigma} =
D_{\sigma(a)} $.

On the other hand, the quotient of $( D') ^{\sigma} $ by the action of the group $G_a$  has genus equal to
the dimension of the space of invariants 
$ dim (H^0 (\Omega^1_{  ( D') ^{\sigma}})^{G_a}  )$, but this dimension is the same $g' =  dim (H^0
(\Omega^1_{  D'})^{G_a}  )$. Hence the action of $G_a$ on $( D') ^{\sigma} $ is also free (by Hurwitz'
formula), and we have shown that $ S^{\sigma}$ is a surface whose moduli point is in
$ \frak M_{\sigma(a)}$.

\qed

\section{Proof of the main theorems}

{\em Prooof of theorem 1.}

Given $ a \in \bar{\Q}$, consider a connected component $\frak N_a$ of $\frak M_a$. Our  main theorem
1 follows easily from the following

{\bf Main Claim:} if $\frak N_a  = \sigma (\frak N_a )$, then necessarily $a  = \sigma (a)$.

\Proof Denote $\sigma (a)$ by $b$. We know that $ C_a \cong D_a / H_a$,  that $ C_{b} \cong D_{b} /
H_{b}$, and we have already shown that $\sigma$ yields an isomorphism $G_a \cong G_{b}$. The
assumption that $\frak N_a  = \sigma (\frak N_a )$ implies, by theorem (3.3) of \cite{cat03},
 the condition that the pairs $(D_a, G_a)$ and $(D_{b}, G_{b})$ are isomorphic as complex triangle curves.
It suffices therefore to show that under this isomorphism the subgroup $H_a$ corresponds to the subgroup
$ H_{b}$, because then we conclude that the curves $C_a$, $C_{b}$ are not only Galois conjugate, but also
isomorphic.

And then  by lemma 1 we conclude that $ a = b$.

To show that $H_a$ corresponds to the subgroup $ H_{b}$,  let $K$ be the Galois closure of the field $L$
(= splitting field of the field extension $\Q \subset L$), and view $L$ as embedded in $\C$ under the
isomorphism sending $x \mapsto a$.

Consider the curve $\hat{C}_x$  obtained from $C_x$ by scalar extension $\hat{C}_x : = C_x \otimes_L K$.
Let also $\hat{F}_x : = F_x \otimes_L K$ the corresponding Belyi function with values in $\PP^1_K$.

Apply now the effective construction of the normal closure of  section \ref{nc}, and, taking a connected
component of  
$ (\hat{C}_x \times _{\PP^1_K} \ldots \times _{\PP^1_K} \hat{C}_x) \setminus \Delta $  we obtain  a
curve $D_x$ defined over $K$.

Note that $D_x$ is not geometrically irreducible, but once we tensor with $\C$ it splits into several
components which are Galois conjugate and which are isomorphic to the conjugates of $D_a$.

Apply now the Galois automorphism $\sigma$ to the triple $D_a \to C_a \to \PP^1$. Since the triple is
induced by the triple  $D_x \to C_x \to \PP^1_K$ by taking a tensor product $ \otimes_K \C$ via the
embedding sending $ x \mapsto a$, the morphisms are induced by the composition of the inclusion $D_x
\subset (C_x)^d$ with the coordinate projections, respectively buy the fibre product equation, it follows
from proposition \ref{Belyi} that $\sigma$ carries the triple $D_a \to C_a \to \PP^1$ to the triple $D_b \to
C_b \to \PP^1$.

\qed

If we want to interpret our argument in terms of Grothendieck's \'etale fundamental group, we  define
$C^0_x : = F_x^{-1} (\PP^1 \setminus \{0,1,\infty \})$, and accordingly  $\hat{C}^0_x$ and $D^0_x$.

There are the following  exact sequences for the  Grothendieck \'etale fundamental group (compare
Theorem 6.1 of \cite{sga1}):

$$ 1 \ra \pi_1^{alg} (D_a) \ra  \pi_1^{alg} (D_x) \ra Gal (\bar{\Q}/ K) \ra 1$$
$$ 1 \ra \pi_1^{alg} (C_a) \ra  \pi_1^{alg} (C_x) \ra Gal (\bar{\Q}/ K)  \ra 1$$
$$ 1 \ra \pi_1^{alg} (\PP^1_{\C}  \setminus \{0,1,\infty \}) \ra  
\pi_1^{alg} (\PP^1_{K}  \setminus \{0,1,\infty \}) \ra Gal (\bar{\Q}/ K) \ra 1$$

where $H_a$ and $G_a$ are the respective factor groups for the inclusions of the left hand sides,
corresponding to the first and second sequence, and  corresponding to the first and third sequence.

On the other hand, we also have the exact sequence

$$ 1 \ra \pi_1^{alg} (\PP^1_{\C}  \setminus \{0,1,\infty \}) \ra  
\pi_1^{alg} (\PP^1_{\Q}  \setminus \{0,1,\infty \}) \ra Gal (\bar{\Q}/ \Q) \ra 1.$$

The finite quotient $G_a$ of $\pi_1^{alg} (\PP^1_{\C}  \setminus \{0,1,\infty \})$ (defined over $K$) is
sent by $\sigma \in Gal (\bar{\Q}/ \Q)$ to another quotient, corresponding to $D_{\sigma (a)}$, and the
subgroup $H_a$, yielding the quotient
$C_a$, is sent to the subgroup  $H_{\sigma (a)}$.

\medskip {\em Proof of theorem 2.}

Assume that $\sigma (a) = b$, and that $\sigma $ is neither complex conjugation nor the identity.
 Consider
a surface $X_a$ with $[X_a]  \in \mathfrak M_a$ and its Galois conjugate 
$ (X_a)^{\sigma } = X_{\sigma (a)} = X_b .$

Denote as before by $\mathfrak N_a$ the connected component of the moduli
 space of surfaces of general type containing $[X_a ]$.

Assume now that the $X_a$ and $X_b$ have isomorphic fundamental groups. Since obviously the two
surfaces have the same Euler number $ e(X_a) = e ( (X_a)^{\sigma } )$ (since they have the same Hodge
numbers $h^i (\Omega^j)$) we can apply again theorem 3.3. of \cite{cat03} and there are only two
possibilities : either $[X_a]$ and $[X_b]$ belong to the same connected component  $\mathfrak N_a$ of
the moduli space, or $[X_b]$ belongs to the complex conjugate  component $\mathfrak N_{\bar a}$ of the
connected component  $\mathfrak N_a$.

In the first case the main claim says that $\sigma (a) = a$, in the second case 
 it  says that $\sigma (a) = \bar a$.

Since however $\sigma $ is neither complex conjugation nor the identity, we find an $a \in \bar{\Q}$ such
that $b : = \sigma (a) \neq a $ and $ b : = \sigma (a) \neq  \bar a$.

Therefore the corresponding surfaces $X_a$ and $(X_a)^{\sigma } $  have nonisomorphic fundamental
groups.

\qed

\begin{oss} We observe that $X_a$ and $(X_a)^{\sigma}$ have isomorphic
 Grothendieck \'etale fundamental groups. In particular, the profinite completion of $\pi_i(X_a)$ and
$\pi_1((X_a)^{\sigma})$ are isomorphic.
\end{oss}

It is not so easy to calculate explicitly the fundamental groups of the
 surfaces constructed above, since one has to explicitly  calculate the monodromy of the Belyi function  of
the very special hyperelliptic curves.

Therefore we give in the next section an explicit example of two  rigid surfaces with
non isomorphic fundamental groups which are Galois conjugate.

\section{An explicit example}

In this section we provide, as we already mentioned, an explicit example of two  surfaces with
non isomorphic fundamental groups which are conjugate under the absolute Galois group, hence with
isomorphic profinite completions of their  respective fundamental groups. These surfaces are rigid.

We consider (see \cite{almeria} for an elementary treatment of what follows) polynomials 
 with only two critical values:
$\{0, 1 \}$.

Let $P \in \mathbb{C}[z]$ be a polynomial with critical values
$\{0,1\}$. 

In order not to have infinitely many polynomials with the same branching behaviour, one considers {\em
normalized polynomials}
 $P(z):= z^n + a_{n-2}z^{n-2} + \ldots a_0$. The condition that $P$ has only $\{0, 1\}$ as critical values,
implies, as we shall briefly recall, that $P$ has coefficients in $\bar{\mathbb{Q}}$. Denote by $K$ the
number field generated by the coefficients of $P$. 

Fix the types the types $(m_1, \ldots, m_r)$ and $(n_1, \ldots, n_s)$ of the cycle
decompositions of the respective local monodromies around $0$ and $1$: we can write our polynomial $P$
in two ways, namely as:
$$
P(z) = \prod_{i=1}^r (z - \beta_i)^{m_i},
$$ and 
$$P(z)  = 1 + \prod_{k=1}^s (z - \gamma_k)^{n_k}.$$

We have the equations $F_1 = \sum m_i \beta_i = 0$ and $F_2=\sum n_k \gamma_k = 0$ (since $P$ is
normalized). Moreover, $m_1 + \ldots + m_r = n_1 + \ldots + n_s = n = degP$ and therefore, since
$\sum_j (m_j-1) + \sum_i(n_i -1) = n-1$, we get $r+s = n+1$.

Since we have
$\prod_{i=1}^r (z -
\beta_i)^{m_i} = 1+ \prod_{k=1}^s (z - \gamma_k)^{n_k}$, comparing coefficients we obtain further
$n-1$ polynomial equations with integer coefficients in the variables $\beta_i$, $\gamma_k$, which we
denote by $F_3=0, \ldots, F_{n+1}=0$. Let $\mathbb{V}(n;(m_1,\ldots,m_n),(n_1,\ldots,n_s))$ be the
algebraic set in affine $(n+1)$-space defined by the equations $F_1=0, \ldots,F_{n+1}=0$. Mapping a
point of this algebraic set to the vector $(a_0,\ldots,a_{n-2})$ of coefficients of the corresponding
polynomial $P$ we obtain a set 
$$\mathbb{W}(n;(m_1,\ldots,m_n),(n_1,\ldots,n_s))$$ (by elimination of variables) in affine $(n-1)$
space. Both these are finite algebraic sets  defined over $\Q$ since by Riemann's existence theorem  they
are either empty or have dimension $0$. 

Observe also that the equivalence classes of monodromies $ \mu \colon \pi_1 (\PP^1 \setminus \{ 0.1.\infty \})
\ra \FS_n$ correspond to the orbits of the group of n-th roots of $1$
(we refer to \cite{almeria} for more details).

We recall also the following
\begin{df}
A smooth algebraic curve $C$ is called a {\em triangle curve} iff there is a finite group $G$ acting
effectively on $C$ with the property  that $C/G \cong \PP^1$ in such a way  that 
$f \colon C \rightarrow \PP^1$ is ramified only in $\{0,1, \infty \}$.
\end{df}

\begin{ex}\label{ex} We calculate (here and in the following, either by a MAGMA routine, or, sometimes,
more painfully by direct
calculation) that
$\mathbb{W}:=\mathbb{W}(7;(2,2,1,1,1);(3,2,2))$ is irreducible over
$\mathbb{Q}$. This implies that
$Gal(\bar{\mathbb{Q}}/\mathbb{Q})$ acts transitively on
$\mathbb{W}$. Looking at the possible monodromies,  one sees that there are exactly two real non
equivalent polynomials. Observe also that the equivalence classes of monodromies $\mu
\colon \pi_1(\PP^1 \setminus \{0,1, \infty \}) \rightarrow \mathfrak(S)_n$ correspond to the orbits of the
group $\mu_n$ of $n$-th roots of $1$ (we refer to
\cite{almeria} for more details). In both cases, which will be explicitly described later on, the two
permutations, of types
$(2,2)$ and
$(3,2,2)$, are seen to generate
$\mathfrak{A}_7$ and the respective normal closures of the two polynomial maps 
are easily seen to give (we use here the fact that the automorphism group of $\mathfrak{A}_7$ is
$\mathfrak{S}_7$) nonequivalent  triangle curves $C_1$, $C_2$. 

By Hurwitz's formula, we see that $g(C_i) = \frac{|\mathfrak{A}_7|}{2}(1 - \frac{1}{2} - \frac{1}{6}
-\frac{1}{7}) + 1 = 241$.

We remark that $\mathfrak{A}_7$ admits generators $a_1, a_2$ of order $5$ such that their product has
order $5$, hence we get a triangle curve
$C$ (of genus $505$). 
\end{ex}

Consider the triangle curve $C$ given by a {\bf spherical system of generators of type $(5,5,5)$ of
$\mathfrak{A}_7$}, i.e., generators $a_1, a_2, a_3$ of $\mathfrak{A}_7$ such that $a_1 \cdot a_2 \cdot
a_3 = 1$. 

\begin{df}
Let $(a_1,a_2,a_3)$ and $(b_1,b_2,b_3)$ be two spherical systems of generators of a finite group $G$ of
the same unordered type, i.e., $\{ord(a_1), ord(a_2), ord(a_3) \} = \{ord(b_1), ord(b_2), ord(b_3) \}$.
Then $(a_1,a_2,a_3)$ and $(b_1,b_2,b_3)$ are called {\bf Hurwitz equivalent} iff they are equivalent
under the equivalence relation generated by
$$
(a_1,a_2,a_3) \equiv (a_2, a_2^{-1}a_1a_2, a_3),
$$
$$
(a_1,a_2,a_3) \equiv (a_1,a_3, a_3^{-1}a_2a_3).
$$
\end{df}

It is well known that two such triangle curves are isomorphic,
compatibly with the action of the group $G$, if and only if the two spherical
systems of generators are {\em Hurwitz equivalent}.  

\begin{oss}
An easy MAGMA routine shows that there is exactly one
Hurwitz equivalence  class of triangle curves given by 
a {\em spherical system of generators of type $(5,5,5)$ of
$\mathfrak{A}_7$}. In other words, if $D_1$ and $D_2$ are two triangle curves given by spherical systems
of generators of type
$(5,5,5)$ of
$\mathfrak{A}_7$, then $D_1$ and $D_2$ are not only isomorphic as algebraic curves, but they have the
same action of $G$.
\end{oss}
 
Let $C$ be as above the triangle curve given by a(ny) spherical systems
of generators of type
$(5,5,5)$ of
$\mathfrak{A}_7$, and consider the two triangle curves $C_1$ and $C_2$ as in example \ref{ex}.
Clearly $\mathfrak{A}_7$ acts freely on $C_1 \times C$ as well as on
$C_2 \times C$ and we obtain two non isomorphic so-called {\em Beauville surfaces} 
$S_1 : = (C_1 \times C) / G$, $S_2 : = (C_2 \times C) / G$.

 Obviously, these two surfaces have the same topological Euler characteristic. If they had isomorphic
fundamental groups, by theorem 3.3 of \cite{cat03}, $S_2$ would be the complex conjugate surface of
$S_1$. In particular,  $C_1$ would be the complex conjugate triangle curve
of  $C_2$: but this is absurd since  we shall show that both $C_1$ and
$C_2$ are real triangle curves. 

\begin{prop} There is a field automorphism $\sigma \in Gal(\bar{\Q}/\Q)$ such that $S_2 =
(S_1)^{\sigma}$.
\end{prop}

\Proof We know that $(S_1)^{\sigma} = ((C_1)^{\sigma} \times (C)^{\sigma}) /G$.  Since there is only one
Hurwitz class of triangle curves given by a spherical system  of generators of type
$(5,5,5)$ of $\mathfrak{A}_7$, we have $(C)^{\sigma} \cong C$ (with the same action of $G$).

\qed

We determine now explicitly the respective fundamental groups of $S_1$ and $S_2$.

In general, let $(a_1, \ldots, a_n)$ and $(b_1, \ldots, b_m)$ be two sets of spherical generators of a finite
group $G$ of respective order types $r:=(r_1, \ldots, r_n)$,
$s:=(s_1, \ldots, s_m)$.  We denote the corresponding `polygonal' curves by $D_1$, resp. $D_2$.

Assume now that the diagonal action of $G$ on $D_1 \times D_2$ is free. We get then the smooth surface
$S:= ( D_1 \times D_2 ) /G$, isogenous to a product.

Denote by $T_r := T(r_1, \ldots, r_n)$ the {\em polygonal group} 
$$
\langle x_1, \ldots, x_{n-1} | x_1^ {r_1} =
\ldots = x_{n-1}^ {r_{n-1}} = (x_1 x_2 \ldots  x_{n-1})^ {r_n} = 1\rangle.
$$ 

We have the exact sequence (cf. \cite{isogenous} cor. 4.7)
$$ 1 \rightarrow \pi_1 \times \pi_2 \rightarrow T_r \times T_s
\rightarrow G \times G \rightarrow 1,
$$ where $\pi_i := \pi_1(D_i)$.

Let $\Delta_G$ be the diagonal in $G \times G$ and let $H$ be the inverse image of $\Delta_G$ under
$\Phi: T_r \times T_s
\rightarrow G \times G$.  We get the exact sequence 
$$ 1 \rightarrow \pi_1 \times \pi_2 \rightarrow H
\rightarrow G \cong \Delta_G \rightarrow 1.
$$ 

\begin{oss}
$\pi_1(S) \cong H$ (cf. \cite{isogenous} cor. 4.7).
\end{oss}

We choose now an arbitrary spherical system of generators of type $(5,5,5)$ of $\mathfrak{A}_7$, 
for instance
$((1,7,6,5,4), (1,3,2,6,7), (2,3,4,5,6))$. Note that we use here MAGMA's notation, where permutations act
on the right (i.e., $ab$ sends $x$ to $(xa)b$).

\begin{figure}[htbp]
\begin{center}
\input{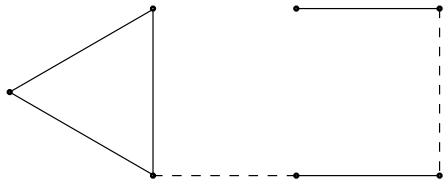}
\end{center}
\caption{Monodromy corresponding to (1)}
\label{figura1}
\end{figure}

A MAGMA routine shows that 
\begin{equation}
((1,2)(3,4), (1,5,7)(2,3)(4,6), (1,7,5,2,4,6,3))
\end{equation}
 and 
\begin{equation}
((1,2)(3,4), (1,7,4)(2,5)(3,6), (1,3,6,4,7,2,5))
\end{equation}
are two representatives of spherical generators of type $(2,6,7)$ yielding
two non isomorphic  triangle curves $C_1$ and $C_2$, each of which
is isomorphic to its complex conjugate. In fact, an alternative
 direct argument is as follows: first of all $C_i$ is isomorphic to its complex
conjugate triangle curve since, for an appropriate choice of the real base point,
complex conjugation sends $ a \mapsto a^{-1}, b \mapsto b^{-1}$
and then one sees that the two corresponding monodromies
are permutation equivalent (see  Figure \ref{figura1} and  Figure \ref{figura2}).

Moreover, since $Aut(\mathfrak{A}_7) = 
\mathfrak{S}_7$, if the two triangle curves were isomorphic, then the two monodromies were conjugate
in $\mathfrak{S}_7$. That this is not the case is seen again by the following pictures. 

\begin{figure}[htbp]
\begin{center}
\input{figura1.pstex_t}
\end{center}
\caption{Monodromy corresponding to (2)}
\label{figura2}
\end{figure}

The two corresponding
homomorphisms
$\Phi_1 :  T_{(2,6,7)} \times T_{(5,5,5)}
\rightarrow \mathfrak{A}_7 \times \mathfrak{A}_7$ and $\Phi_2 :  T_{(2,6,7)} \times T_{(5,5,5)}
\rightarrow \mathfrak{A}_7 \times \mathfrak{A}_7$ give two exact sequences

$$ 1 \rightarrow \pi_1(C_1) \times \pi_1(C) \rightarrow T_{(2,6,7)} \times T_{(5,5,5)}
\rightarrow \mathfrak{A}_7 \times \mathfrak{A}_7 \rightarrow 1,
$$ and
$$ 1 \rightarrow \pi_1(C_2) \times \pi_1(C) \rightarrow T_{(2,6,7)} \times T_{(5,5,5)}
\rightarrow \mathfrak{A}_7 \times \mathfrak{A}_7 \rightarrow 1,
$$

yielding two non isomorphic fundamental groups $\pi_1(S_1) = \Phi_1 ^{-1} (\Delta_{\mathfrak{A}_7})$
and $\pi_1(S_2) = \Phi_2 ^{-1} (\Delta_{\mathfrak{A}_7})$ fitting both in an exact sequence of type
$$ 1 \rightarrow \Pi_{241} \times \Pi_{505}  \rightarrow \pi_1(S_j) 
\rightarrow
\Delta_{\mathfrak{A}_7} \cong \mathfrak{A}_7 \rightarrow 1,
$$

where $\Pi_{241} \cong \pi_1(C_1) \cong \pi_1(C_2)$, $\Pi_{505} = \pi_1(C)$.

\begin{oss} 
1) Using the same trick that we used for our main theorems, namely, using a surjection
of a group $\Pi_g \ra \mathfrak{A}_7$ we get infinitely many examples of pairs of
fundamental groups which are nonisomorphic, but which have 
isomorphic profinite completions. Each pair fits into an exact sequence
$$ 1 \rightarrow \Pi_{241} \times \Pi_{g'}  \rightarrow \pi_1(S_j) 
\rightarrow
 \mathfrak{A}_7 \rightarrow 1.
$$

2)  Many more explicit examples as the one above (but with cokernel group different from 
${\mathfrak{A}_7}$) can be obtained
using polynomials with two critical values.

A construction of polynomials with two critical values having a
very large Galois orbit was proposed to us by D. van Straten.
\end{oss}

\begin{oss}
In the sequel to this work we plan to show that the absolute Galois group in fact
acts faithfully also on the set of Beauville surfaces.
\end{oss}

{\bf Acknowledgements.} The research of the authors was performed in the realm of the
Forschergruppe 'Classification of algebraic surfaces and compact complex manifolds' of the D.F.G.. We
would like to thank Ravi Vakil for his interest in our work and for pointing
out a minor error in a first version of this note.

\noindent {\bf Authors' addresses:}

\noindent Ingrid Bauer, Fabrizio Catanese\\ Lehrstuhl Mathematik VIII, Universit\"at Bayreuth\\
  D-95440 Bayreuth, Germany\\
\noindent (email: Ingrid.Bauer@uni-bayreuth.de, Fabrizio.Catanese@uni-bayreuth.de)\\

\noindent Fritz Grunewald\\ Mathematisches Institut der Heinrich-Heine Universit\"at D\"usseldorf \\
 Universit\"atsstrasse \\ D\"usseldorf \\ 
\noindent (email: grunewald@math.uni-duesseldorf.de)

\end{document}